
\magnification=\magstep1
\parskip 1pt plus2pt minus2pt
\hoffset=0.1 truecm
\voffset=-0.25 truecm
\hsize=15.5 truecm
\vsize=24.5 truecm
%
%
\def\bbbq{{\mathchoice
{\setbox0=\hbox {$\displaystyle\rm Q$}\hbox
{\raise0.15\ht0\hbox to0pt{\kern0.4\wd0\vrule height0.8\ht0\hss}\box0}}
{\setbox0=\hbox {$\textstyle\rm Q$}\hbox
{\raise0.15\ht0\hbox to0pt{\kern0.4\wd0\vrule height0.8\ht0\hss}\box0}}
{\setbox0=\hbox {$\scriptstyle\rm Q$}\hbox
{\raise0.15\ht0\hbox to0pt{\kern0.4\wd0\vrule height0.7\ht0\hss}\box0}}
{\setbox0=\hbox {$\scriptscriptstyle\rm Q$}\hbox
{\raise0.15\ht0\hbox to0pt{\kern0.4\wd0\vrule height0.7\ht0\hss}\box0}}
}}


\def\bbbc{{\mathchoice
{\setbox0=\hbox {$\displaystyle\rm C$}\hbox
{\hbox to0pt{\kern0.4\wd0\vrule height0.9\ht0\hss}\box0}}
{\setbox0=\hbox {$\textstyle\rm C$}\hbox
{\hbox to0pt{\kern0.4\wd0\vrule height0.9\ht0\hss}\box0}}
{\setbox0=\hbox {$\scriptstyle\rm C$}\hbox
{\hbox to0pt{\kern0.4\wd0\vrule height0.9\ht0\hss}\box0}}
{\setbox0=\hbox {$\scriptscriptstyle\rm C$}\hbox
{\hbox to0pt{\kern0.4\wd0\vrule height0.9\ht0\hss}\box0}}
}}

\font\fivesans=cmss10 at 5pt
\font\sevensans=cmss10 at 7pt
\font\tensans=cmss10
\newfam\sansfam
\textfont\sansfam=\tensans\scriptfont\sansfam=\sevensans
\scriptscriptfont\sansfam=\fivesans
\def\sans{\fam\sansfam\tensans}
\def\bbbz{{\mathchoice {\hbox{$\sans\textstyle Z\kern-0.4em Z$}}
{\hbox{$\sans\textstyle Z\kern-0.4em Z$}}
{\hbox{$\sans\scriptstyle Z\kern-0.3em Z$}}
{\hbox{$\sans\scriptscriptstyle Z\kern-0.2em Z$}}}}


\font \bigbf=cmbx10 scaled \magstep2

\def\slash#1{#1\kern-0.65em /}
\def\dirac{{\raise0.09em\hbox{/}}\kern-0.58em\partial}
\def\Dirac{{\raise0.09em\hbox{/}}\kern-0.69em D}





%
\vglue 1.5cm

\centerline {\bigbf On Curvature in Noncommutative Geometry}
\vskip 1.5cm

\centerline {\bf M. Dubois-Violette, \ J. Madore, \ T. Masson}
\medskip
\centerline {\it Laboratoire de Physique Th\'eorique et Hautes
Energies\footnote{*}{\it Laboratoire associ\'e au CNRS.}}
\centerline {\it Universit\'e de Paris-Sud, B\^at. 211,  \ F-91405 Orsay}
\vskip 1cm

\centerline {\bf J. Mourad}
\medskip
\centerline {\it  D\'epartement de Physique, 49 av. des Genottes,}
\centerline {\it Universit\'e de Cergy Pontoise, BP 8428,
                 F-95806 Cergy Pontoise}
\vskip 2cm

\noindent
{\bf Abstract:} \ A general definition of a bimodule connection in
noncommutative geometry has been recently proposed. For a given algebra
this definition is compared with the ordinary definition of a connection
on a left module over the associated enveloping algebra. The
corresponding curvatures are also compared.

\vfill
\noindent
LPTHE Orsay 95/63
\medskip
\noindent
October, 1995
\bigskip
\eject

\beginsection 1 Introduction and motivation

Recently a general definition has been given (Mourad 1995,
Dubois-Violette {\it et al.} 1995) of a linear connection in the context
of noncommutative geometry which makes essential use of the full
bimodule structure of the differential forms. A preliminary version of
the curvature of the connection was given (Madore {\it et al.} 1995)
which had the drawback of not being in general a linear map with respect
to the right-module structure. It is in fact analogous to the curvature
which is implicitly used by those authors (Chamseddine {\it et al.}
1993, Sitarz 1994, Klim\v c\' ik {\it et al.} 1994, Landi {\it et al.}
1994,) who define a linear connection using the formula for a covariant
derivative on an arbitrary left (or right) module (Karoubi 1981, Connes
1986). Our purpose here is to present a modified definition of curvature
which is bilinear. Let ${\cal A}$ be a general associative algebra (with
unit element). This is what replaces in noncommutative geometry the
algebra of smooth functions on a smooth (compact) manifold which is used
in ordinary differential geometry.  By `bilinear' we mean, here and in
what follows, bilinear with respect to ${\cal A}$. In fact we shall
present two definitions of curvature.  The first is valid in all
generality and reduces to the ordinary definition of curvature in the
commutative case.  The second definition seems to be better adapted to
`extreme' noncommutative cases, such as the one considered in Section~5.

The definition of a connection as a covariant derivative was given an
algebraic form in the Tata lectures by Koszul (1960) and generalized to
noncommutative geometry by Karoubi (1981) and Connes (1986, 1994). We
shall often use here the expressions `connection' and `covariant
derivative' synonymously.  In fact we shall distinguish three different
types of connections. A `left ${\cal A}$-connection' is a connection on
a left ${\cal A}$-module; it satisfies a left Leibniz rule.  A `bimodule
${\cal A}$-connection' is a connection on a general bimodule ${\cal M}$
which satisfies a left and right Leibniz rule. In the particular case
where ${\cal M}$ is the module of 1-forms we shall speak of a `linear
connection'.  The precise definitions are given below.  A bimodule over
an algebra ${\cal A}$ is also a left module over the tensor product
${\cal A}^e = {\cal A} \otimes_\bbbc {\cal A}^{\rm op}$ of the algebra
with its `opposite'. So a bimodule can have a bimodule
${\cal A}$-connection as well as a left ${\cal A}^e$-connection. These
two definitions are compared in Section~2. In Section~3 we discuss the
curvature of a bimodule connection.  In Section~4 we consider an algebra
of forms based on derivations and we compare the left connections with
the linear connections. We show that in a sense to be made precise the
two definitions yield the same bilinear curvature. That is, the extra
restriction which the bimodule structure seems to place on the linear
connections does not in fact restrict the corresponding curvature. In
Section~5 we consider a more abstract geometry whose differential
calculus is not based on derivations. In Section~6 a possible definition
is given of the curvature of linear connections on braided-commutative
algebras. In Section~7 we examine the (left) projective structure of the
1-forms of the Connes-Lott model.

Let ${\cal A}$ be an arbitrary algebra and $(\Omega^*({\cal A}) ,d)$ a
differential calculus over ${\cal A}$. One defines a left
${\cal A}$-connection on a left ${\cal A}$-module ${\cal H}$ as a
covariant derivative
$$
{\cal H} \buildrel D \over \rightarrow
\Omega^1({\cal A}) \otimes_{\cal A} {\cal H}                      \eqno(1.1)
$$
which satisfies the left Leibniz rule
$$
D (f \psi) =  df \otimes \psi + f D\psi                           \eqno(1.2)
$$
for arbitrary $f \in {\cal A}$. This map has a natural extension
$$
\Omega^*({\cal A}) \otimes_{\cal A} {\cal H}
\buildrel \nabla \over \longrightarrow
\Omega^*({\cal A}) \otimes_{\cal A} {\cal H}                      \eqno(1.3)
$$
given, for $\psi \in {\cal H}$ and $\alpha \in \Omega^n({\cal A})$,
by $\nabla \psi = D \psi$ and
$$
\nabla (\alpha \psi) = d\alpha \otimes \psi +
                (-1)^n \alpha \nabla \psi.
$$
The operator $\nabla^2$ is necessarily left-linear. However when
${\cal H}$ is a bimodule it is not in general right-linear.

A covariant derivative on the module $\Omega^1({\cal A})$ must satisfy
(1.2). But $\Omega^1({\cal A})$ has also a natural structure as a right
${\cal A}$-module and one must be able to write a corresponding right
Leibniz rule in order to construct a bilinear curvature. Quite generally
let ${\cal M}$ be an arbitrary bimodule. Consider a covariant derivative
$$
{\cal M} \buildrel D \over \rightarrow
\Omega^1({\cal A}) \otimes_{\cal A} {\cal M}                       \eqno(1.4)
$$
which satisfies both a left and a right Leibniz rule.  In order to
define a right Leibniz rule which is consistent with the left one, it
was proposed by Mourad (1995), by Dubois-Violette \& Michor (1995) and
by Dubois-Violette \& Masson (1995) to introduce a generalized
permutation
$$
{\cal M} \otimes_{\cal A} \Omega^1({\cal A})
\buildrel \sigma \over \longrightarrow
\Omega^1({\cal A}) \otimes_{\cal A} {\cal M}.                     \eqno(1.5)
$$
The right Leibniz rule is given then as
$$
D(\xi f) = \sigma (\xi \otimes df) + (D\xi) f                     \eqno(1.6)
$$
for arbitrary $f \in {\cal A}$ and $\xi \in {\cal M}$. The purpose of
the map $\sigma$ is to bring the differential to the left while
respecting the order of the factors. It is necessarily bilinear
(Dubois-Violette {\it et al.} 1995). We define a bimodule
${\cal A}$-connection to be the couple $(D, \sigma)$.

If in particular
$$
{\cal M} = \Omega^1({\cal A})                                     \eqno(1.7)
$$
then we shall refer to the bimodule ${\cal A}$-connection as a linear
connection.  Although we shall here be concerned principally with this
case we shall often consider more general situations.  In any case we
shall use the more general notation to be able to distinguish the two
copies of $\Omega^1({\cal A})$ on the right-hand side of (1.4).

Let $\Omega_u^*({\cal A})$ be the universal differential calculus.
Dubois-Violette \& Masson (1995) have shown that given an arbitrary left
connection on a bimodule ${\cal M}$ there always exists a bimodule
homomorphism
$$
{\cal M} \otimes_{\cal A} \Omega_u^1({\cal A})
\buildrel \sigma(D) \over \longrightarrow
\Omega^1({\cal A}) \otimes_{\cal A} {\cal M}
$$
such that
$$
D(\xi f) = \sigma(D) (\xi \otimes d_uf) + (D\xi) f.
$$
The notation $\sigma(D)$ is taken from the definition of the symbol
of a differential operator. The condition (1.6) means then that
$\sigma(D)$ factorizes as a composition of a $\sigma$ as above and
the canonical homomorphism of
${\cal M} \otimes_{\cal A} \Omega_u^1({\cal A})$ onto
${\cal M} \otimes_{\cal A} \Omega^1({\cal A})$.

Using $\sigma$ one can also construct (Mourad 1995) an extension
$$
\Omega^1({\cal A}) \otimes_{\cal A} \Omega^1({\cal A})
\buildrel D \over \longrightarrow
\Omega^1({\cal A}) \otimes_{\cal A} \Omega^1({\cal A})
\otimes_{\cal A} \Omega^1({\cal A})                                \eqno(1.8)
$$
It can also be proved in fact (Bresser {\it et al.} 1995) that this
extension implies the existence of $\sigma$.  The operator $D^2$
is not in general left-linear. However if we define $\pi$ to be the
product in $\Omega^*({\cal A})$ and set $\pi_{12} = \pi \otimes 1$ then
$\pi_{12} D^2$ is left-linear,
$$
\pi_{12} D^2 (f \xi) = f \pi_{12} D^2 \xi,                         \eqno(1.9)
$$
provided the torsion vanishes and the map $\sigma$ satisfies the
condition
$$
\pi \circ (\sigma + 1) = 0.                                       \eqno(1.10)
$$
The map $\nabla$ is related to $D$ on ${\cal H} = \Omega^1({\cal A})$ by
$$
\nabla^2 = \pi_{12} \circ D^2.                                    \eqno(1.11)
$$
The left-hand side of this equation is define for a general
${\cal A}$-connection whereas the right-hand side is defined only in the
case of a linear connection.

The torsion $T$ is defined to be the map
$$
T = d - \pi \circ D                                              \eqno(1.12)
$$
from $\Omega^1$ into $\Omega^2$.  The restriction (1.7) is here
essential.  It follows from the condition (1.10) that $T$ is bilinear.
A metric can be defined and it can be required to be symmetric using the
map $\sigma$.  The standard condition that the covariant derivative be
metric-compatible can be also carried over to the noncommutative case.
For more details we refer, for example, to Madore {\it et al.} (1995).

\beginsection 2 The bimodule structure

For any algebra ${\cal A}$ the enveloping algebra ${\cal A}^e$ is
defined to be
$$
{\cal A}^e = {\cal A} \otimes_\bbbc {\cal A}^{\rm op}.
$$
A bimodule ${\cal M}$ can also be considered then as a left
${\cal A}^e$-module.  The differential calculus $\Omega^*({\cal A})$ has
a natural extension to a differential calculus $\Omega^*({\cal A}^e)$
given by
$$
\Omega^*({\cal A}^e) =
\Omega^*({\cal A}) \otimes \Omega^*({\cal A}^{\rm op}) =
\big(\Omega^*({\cal A})\big)^e                                    \eqno(2.1)
$$
with $d (a \otimes b) = da \otimes b + a \otimes db$.  This is not the
only choice. For example if $\Omega^*({\cal A})$ were the universal
calculus over ${\cal A}$ then $\Omega^*({\cal A}^e)$ would not be
equal to the universal calculus over ${\cal A}^e$. Suppose that
${\cal M}$ has a left ${\cal A}^e$-connection
$$
{\cal M} \buildrel  D^e \over \longrightarrow
\Omega^1({\cal A}^e) \otimes_{{\cal A}^e} {\cal M}.               \eqno(2.2)
$$
{}From the equality
$$
\Omega^1({\cal A}^e) =
\big(\Omega^1({\cal A}) \otimes_\bbbc {\cal A}^{\rm op}\big) \oplus
\big({\cal A} \otimes_\bbbc \Omega^1({\cal A}^{\rm op})\big).     \eqno(2.3)
$$
and using the identification
$$
\big({\cal A} \otimes_\bbbc \Omega^1({\cal A}^{\rm op})\big)
\otimes_{{\cal A}^e} {\cal M} \simeq
{\cal M} \otimes_{\cal A} \Omega^1({\cal A})                      \eqno(2.4)
$$
given by
$$
(1 \otimes \xi) \otimes \eta \mapsto \eta \otimes \xi
$$
we find that we have
$$
\Omega^1({\cal A}^e) \otimes_{{\cal A}^e} {\cal M} =
\big(\Omega^1({\cal A}) \otimes_{\cal A} {\cal M}\big) \oplus
\big({\cal M} \otimes_{\cal A} \Omega^1({\cal A})\big).           \eqno(2.5)
$$
The covariant derivative $D^e$ splits then as the sum of two terms
$$
D^e = D_L + D_R.                                                  \eqno(2.6)
$$
{}From the identifications it is obvious that $D_L$ ($D_R$) satisfies a
left (right) Leibniz rule and is right (left) ${\cal A}$-linear. Such
covariant derivatives have been considered by Cuntz \& Quillen (1995),
by Bresser {\it et al.} (1995) and by Dabrowski {\it et al.} (1995).

One can write a (noncommutative) triangular diagram
$$

\matrix{                      & {\cal M } &                     \cr
        \hfill  D_L \swarrow  &           & \searrow D_R \hfill \cr
\Omega^1({\cal A}) \otimes_{\cal A} {\cal M}  &
\buildrel \sigma \over \longleftarrow         &
{\cal M} \otimes_{\cal A} \Omega^1({\cal A})}                     \eqno(2.7)

$$
from which one sees that given an arbitrary bimodule homomorphism (1.5)
and a covariant derivative (2.2) one can construct a covariant
derivative (1.4) by the formula
$$
D = D_L + \sigma \circ D_R                                        \eqno(2.8)
$$
which satisfies both (1.2) and (1.6).

Suppose further that the differential calculus is such that the
differential $d$ of an element $f \in {\cal A}$ is of the form
$$
d f = - [\theta, f],                                              \eqno(2.9)
$$
for some element $\theta \in \Omega^1$. Then obviously
particular choices for $D_L$ and $D_R$ are the expressions
$$
D_L \xi = - \theta \otimes \xi, \qquad
D_R \xi = \xi \otimes \theta.                                     \eqno(2.10)
$$
Let $\tau$ be a bimodule homomorphism from ${\cal M}$ into
$\Omega^1({\cal A}^e) \otimes_{{\cal A}^e} {\cal M}$ and decompose
$$
\tau = \tau_L + \tau_R                                            \eqno(2.11)
$$
according to the decomposition (2.5). The most general $D_L$ and $D_R$
are of the form
$$
D_L \xi = - \theta \otimes \xi + \tau_L(\xi), \qquad
D_R \xi = \xi \otimes \theta + \tau_R(\xi).                       \eqno(2.12)
$$

Using (2.8) we can construct a covariant derivative
$$
D\xi = - \theta \otimes \xi + \sigma (\xi \otimes \theta)         \eqno(2.13)
$$
from (2.10). In Section~4 we shall study a differential calculus for
which this is the only possible $D$.

{}From the Formula~(2.9) we know that there is a bimodule projection of
${\cal A}^e$ onto $\Omega^1({\cal A})$. Suppose that
$\Omega^1({\cal A})$ is a projective bimodule and let $P$ be the
corresponding projector. We can identify then $\Omega^1({\cal A})$ as a
sub-bimodule of the free ${\cal A}^e$-module of rank 1:
$$
\Omega^1({\cal A}) = {\cal A}^e \, P.
$$

A left ${\cal A}^e$-connection on ${\cal A}^e$ as a left
${\cal A}^e$-module is a covariant derivative of the form (1.4) with
${\cal M} = {\cal A}^e$.  The ordinary differential $d^e$ on
${\cal A}^e$,
$$
{\cal A}^e \buildrel d^e \over \rightarrow \Omega^1({\cal A}^e),  \eqno(2.14)
$$
is clearly a covariant derivative in this sense.  The right-hand side
can be written using (2.5) as
$$
\Omega^1({\cal A}^e) =
\big(\Omega^1({\cal A}) \otimes_{\cal A} {\cal A}^e\big) \oplus
\big({\cal A}^e \otimes_{\cal A} \Omega^1({\cal A})\big)           \eqno(2.15)
$$
and so we can split $d^e$ as the sum of two terms $d_L$ and $d_R$.
Let $a \otimes b$ be an element of ${\cal A}^e$. Then we have
$$
d_L (a \otimes b) = - [\theta, a] \otimes b =
- \theta \otimes (a \otimes b) + (a \otimes b) (\theta \otimes 1). \eqno(2.16)
$$
In the first term on the right-hand side the first tensor product is
over the algebra and the second is over the complex numbers; in the
second term the first tensor product is over the complex numbers and the
second is over the algebra.

A general element of $\Omega^1({\cal A})$ can be written as a sum
of elements of the form $\xi = (a \otimes b) P = a P b$. We have
then
$$
d_L \xi = - \theta \otimes \xi + \xi (\theta \otimes 1).
$$
Define $D_L$ by
$$
D_L \xi = (d_L \xi) P.                                             \eqno(2.17)
$$
Then we obtain the first of equations (2.12) with
$$
\tau_L (\xi) = \xi (\theta \otimes P).                             \eqno(2.18)
$$
Here, on the right-hand side, the tensor product is over the
algebra and $\theta \otimes P$ is an element of
$\Omega^1({\cal A}) \otimes \Omega^1({\cal A})$. This is a left
${\cal A}^e$-module. Similarly one can construct a $D_R$ and a $D^e$
by equation (2.6):
$$
D^e \xi = (d^e \xi) P.                                            \eqno(2.19)
$$

In the case of ordinary geometry with ${\cal A}$ equal to the algebra
${\cal C}^\infty (V)$ of smooth functions on a smooth manifold $V$ the
algebra ${\cal A}^e$ is the algebra of smooth functions in two
variables.  If $\Omega^*({\cal A})$ is the algebra of de~Rham
differential forms the only possible $\sigma$ is the permutation and the
left and right Leibniz rules are identical.  In this case $D^e$ cannot
exist.  In fact $D_L$ would satisfy a left Leibniz rule and be left
linear since the left and right multiplication are equal. In general let
${\cal M}$ be the ${\cal A}$-module of smooth sections of a vector
bundle over $V$.  Then ${\cal M}$ is a ${\cal A}^e$-module.  It is
important to notice that although it is projective as an
${\cal A}$-module it is never projective as an ${\cal A}^e$-module since
a projective ${\cal A}^e$-module consists of 2-point functions.

\beginsection 3 Curvature

Consider a covariant derivative (1.4) which satisfies the left Leibniz
rule (1.2). We can define a right-linear curvature by factoring out in
the image of $\nabla^2$ all those elements ($J$ = `junk') which do not
satisfy the desired condition. Define $J$ as the vector space
$$
J = \big\{ \sum_i \big(\nabla^2 (\xi_i f_i) - \nabla^2 (\xi_i) f_i \big)
    \; \big\vert \; \xi_i \in {\cal M}, f_i \in {\cal A} \big\}.   \eqno(3.1)
$$
In fact $J$ is a sub-bimodule of
$\Omega^2({\cal A}) \otimes_{\cal A} {\cal M}$. It is obviously a
left-submodule.  Consider the element
$\alpha = \nabla^2 (\xi g ) - \nabla^2 (\xi) g \in J$ and let
$f \in {\cal A}$. We can write
$$
\alpha f = \big(\nabla^2 (\xi g f) - \nabla^2 (\xi) g f\big)
         - \big(\nabla^2 (\xi g f) - \nabla^2 (\xi g) f\big).
$$
Therefore $\alpha f \in J$ and $J$ is also a right submodule.

Let $p$ be the projection
$$
\Omega^2({\cal A}) \otimes_{\cal A} {\cal M}
\buildrel p \over \longrightarrow
\Omega^2({\cal A}) \otimes_{\cal A} {\cal M}/ J.                   \eqno(3.2)
$$
We shall define the curvature of $D$ as the combined map
$$
{\rm Curv} = - p \circ \nabla^2.                                   \eqno(3.3)
$$
In the case of a linear connection we can write
$$
{\rm Curv} = - p \circ \pi_{12} \circ D^2.
$$
By construction Curv is left and right-linear:
$$
{\rm Curv}(f\xi) = f {\rm Curv}(\xi), \qquad
{\rm Curv}(\xi f) = {\rm Curv}(\xi) f.                             \eqno(3.4)
$$
In the next section we shall present an example which illustrates the
role which the right-Leibniz rule (1.6) plays in this construction.

Consider the covariant derivative (2.2). One can define a bilinear
curvature as the map
$$
{\rm Curv}_L = - \nabla_L^2                                         \eqno(3.5)
$$
from ${\cal M}$ into $\Omega^2({\cal A}) \otimes_{\cal A} {\cal M}$. It is
bilinear because by construction it is trivially right-linear.  In the
case where the differential $d$ is given by (2.9) and $D_L$ is given by
(2.10) we find that ${\rm Curv}_L$ is given by the formula
$$
{\rm Curv}_L(\xi) = (d\theta + \theta^2) \otimes \xi.             \eqno(3.6)
$$
{}From this expression it is obvious that ${\rm Curv}_L$ is right-linear;
it is easy to verify directly that it is also left-linear because of the
fact that the 2-form $d\theta + \theta^2$ commutes with the elements of
${\cal A}$:
$$
[d\theta + \theta^2, f] = d[\theta, f] = - d^2 f = 0.              \eqno(3.7)
$$

The covariant derivative (2.2) has also a bilinear curvature 2-form
$$
{\rm Curv}^e = - \nabla^{e2}                                       \eqno(3.8)
$$
which naturally decomposes into three terms all of which are bilinear.
One of these terms corresponds to the covariant derivative of Section~1
with $\sigma$ set equal to zero.  It takes its values in a space which
can be naturally identified with
$\Omega^2({\cal A}) \otimes_{\cal A} {\cal M}$.  However because the
second action of $D^e$ does not commute with that of $\sigma$, the
corresponding term of ${\rm Curv}^e$ does not necessarily coincide
with the image of ${\rm Curv}$. We shall discuss this in an example in
Section~5.  The curvature of the particular connection (2.19) can be
written in terms of the projector $P$:
$$
{\rm Curv}^e \xi = - \nabla^{e2} \xi
                 = - \xi \big((d^e P) (d^e P) P\big).            \eqno(3.9)
$$

The extension of $D$ to the tensor product of $n$ copies of
$\Omega^1({\cal A})$ defines a covariant derivative on the left module
$$
{\cal H} = \big(\Omega^1({\cal A})\big)^{\otimes n}.             \eqno(3.10)
$$
The curvature is given by (1.11). In the commutative case and, more
generally, in the case of a derivation-based differential calculus this
curvature can be expressed in terms of the curvature of the covariant
derivative (1.4). For a general differential calculus this will not be
the case.

The same remarks can be made concerning the torsion (1.12).
In general let $\pi$ be the product map of
$\big(\Omega^1({\cal A})\big)^{\otimes n}$ into $\Omega^n({\cal A})$.
Then one can also define a module homomorphism
$$
\big(\Omega^1({\cal A})\big)^{\otimes n}
\buildrel T_n \over \longrightarrow \Omega^{n+1}({\cal A})       \eqno(3.11)
$$
given by
$$
T_n=d\pi-\pi \circ D.                                            \eqno(3.12)
$$
These maps are all left-module homomorphisms. If
$\xi \in \Omega^1({\cal A})$ and
$\nu \in \big(\Omega^1({\cal A})\big)^{\otimes n}$ then we have
$$
T_{n+1}(\xi \otimes \nu) = T_1 (\xi) \, \pi (\nu) - \xi \, T_n(\nu)
+ \pi \circ \big( (\sigma+1)\otimes 1\big) \xi \otimes \nabla \nu.\eqno(3.13)
$$
In order for the last term in the previous equation to vanish it is
necessary and sufficient that (1.10) be satisfied.  In this case one
sees by iteration that the $T_n$ can all be expressed in terms of $T_1$
and therefore that all of them are bimodule homomorphisms.

\beginsection 4 Linear connections on matrix geometries

As a first example we present the case of the algebra $M_n$ of $n\times n$
matrices (Dubois-Violette {\it et al.} 1989, 1990) with a differential
calculus based on derivations (Dubois-Violette 1988). Let $\lambda_r$ be
a set of generators of the Lie algebra of the special linear group
$SL_n$. Then the derivations $e_r = {\rm ad}\, \lambda_r$ is a basis for
the derivations of $M_n$ and the dual 1-forms $\theta^r$ commute with
the elements of $M_n$.  The set of 1-forms $\Omega^1(M_n)$ is a free
left (or right) module of rank $n^2-1$.  The natural map $\sigma$ which
we shall use is given (Madore {\it et al.} 1995) by
$$
\sigma (\theta^r \otimes \theta^s) = \theta^s \otimes \theta^r.   \eqno(4.1)
$$
Quite generally for any algebra ${\cal A}$ with a differential calculus
which is based on derivations there is a natural map $\sigma$ given by a
permutation of the arguments in the forms. Let $X$ and $Y$ be
derivations. Then one can define $\sigma$ by
$$
\sigma (\xi \otimes \eta) (X, Y) = \xi \otimes \eta (Y, X).
$$

A general left $M_n$-connection can be defined by the covariant
derivative
$$
D\theta^r = -  \omega^r{}_{st} \,\theta^s \otimes \theta^t        \eqno(4.2)
$$
with $\omega^r{}_{st}$ an arbitrary element of $M_n$ for each value of
the indices $r,s,t$.  We write
$$
\omega^r{}_{st} = \Gamma^r_{st} + J^r{}_{st}                      \eqno(4.3)
$$
where the $\Gamma^r_{st}$ are proportional to the identity in $M_n$
and the $J^r{}_{st}$ are trace-free.  If we require that the torsion
vanish then we have (Madore {\it et al.} 1995)
$$
\Gamma^r_{[st]} =  C^r{}_{st}                                     \eqno(4.4)
$$
where the $C^r{}_{st}$ are $SL_n$ structure constants.

If we impose the right Leibniz rule we find that
$$
0 = D ([f,\theta^r]) = [f, D\theta^r]
  = - [f, J^r{}_{st}]\,\theta^s \otimes \theta^t                  \eqno(4.5)
$$
for arbitrary $f \in M_n$ and so we see that if the connection is a
linear connection then
$$
J^r{}_{st} = 0.                                                   \eqno(4.6)
$$

Consider now the curvature of the left $M_n$-connection and write
$$
\nabla^2 \theta^r =
         - \Omega^r{}_{stu} \theta^t \theta^u \otimes \theta^s.   \eqno(4.7)
$$
Then since the elements of the algebra commute with the generators
$\theta^r$ we have
$$
\nabla^2 (\theta^r f) - (\nabla^2 \theta^r) f
= \nabla^2 (f \theta^r) - (\nabla^2 \theta^r) f
= - [f, \Omega^r{}_{stu}] \theta^t \theta^u \otimes \theta^s.    \eqno(4.8)
$$
Since $f$ is arbitrary it follows then that we have
$$
{\rm Curv}(\theta^r) = {1\over 2} R^r{}_{stu}
                            \theta^t \theta^u \otimes \theta^s   \eqno(4.9)
$$
where the $R^r{}_{stu}$ are defined uniquely in terms of the
$\Gamma^r_{st}$:
$$
R^r{}_{stu} = \Gamma^r_{tp} \Gamma^p_{us}
      - \Gamma^r_{up} \Gamma^p_{ts} - \Gamma^r_{ps} C^p{}_{tu}.  \eqno(4.10)
$$
That is, $R^r{}_{stu}$ does not depend on $J^r{}_{st}$.

We conclude then that even had we not required the right Leibniz rule
and had admitted an extra term of the form $J^r{}_{st}$ in the
expression for the covariant derivative then we would find that the
curvature map Curv would remain unchanged. The extra possible terms are
eliminated under the projection $p$ of (3.2).

There is a covariant derivative which is of the form (2.8) with $D_L$
and $D_R$ given by (2.10). For this covariant derivative one has
$$
\omega^r{}_{st} \equiv 0.                                        \eqno(4.11)
$$
This covariant derivative has obviously vanishing curvature but it is
not torsion-free.  If we use the ambiguity (2.11) we can write any
covariant derivative (4.2) in the form (2.8).

The generators $\theta^r$ are no longer independent if one considers the
bimodule structure.  In fact one finds that
$$
\theta^r = - C^r{}_{st} \lambda^s \theta \lambda^t, \qquad
\theta = - \lambda_r \theta^r                                     \eqno(4.12)
$$
and as a bimodule $\Omega^1(M_n)$ is generated by $\theta$ alone. For
dimensional reasons $\Omega^1(M_n)$ cannot be of rank one.  In fact the
free $M_n$-bimodule of rank one is of dimension $n^4$ and the dimension
of $\Omega^1(M_n)$ is equal to $(n^2-1)n^2 < n^4$. With the normalization
which we have used for the generators $\lambda_r$ the element
$$
\zeta = {1\over n^2} 1 \otimes 1 - {1\over n} \lambda_r \otimes \lambda^r
$$
is a projector in $M_n \otimes M_n$ which commutes with the elements of
$M_n$. This can be written as
$$
d(M_n) \zeta = 0.
$$
We have the direct-sum decomposition
$$
M_n \otimes M_n = \Omega^1(M_n) \oplus M_n \,\zeta.              \eqno(4.13)
$$
As in Section~1 one can define $M_n^e = M_n \otimes_\bbbc M_n^{\rm op}$.
The prescription (2.19) with
$$
P = 1 \otimes 1 - \zeta
$$
yields then a covariant derivative of the form (4.2) whose curvature
vanishes.

\beginsection 5 Linear connections on the Connes-Lott model

Consider the algebra $M_3$ with the grading defined by the decomposition
$\bbbc^3 = \bbbc^2 \oplus \bbbc$. Define (Connes \& Lott 1990, 1992)
$\Omega^0(M_3^+) = M_3^+ = M_2 \times M_1$, $\Omega^1(M_3^+) = M_3^-$,
$\Omega^2(M_3^+) = M_1$ and $\Omega^p(M_3^+) = 0$ for $p\geq 3$. A
differential $d$ can be defined by (Connes 1986, 1990)
$$
df = - [\theta, f],                                              \eqno(5.1)
$$
where $\theta \in \Omega^1(M_3^+)$.

The vector space of 1-forms is of dimension 4 over the complex numbers.
The dimension of $\Omega^1(M_3^+) \otimes_\bbbc\Omega^1(M_3^+)$ is equal
to 16 but the dimension of the tensor product
$\Omega^1(M_3^+) \otimes_{M_3^+}\Omega^1(M_3^+)$ is equal to 5 and we
can make the identification
$$
\Omega^1(M_3^+) \otimes_{M_3^+}\Omega^1(M_3^+) = M_3^+.           \eqno(5.2)
$$

To define a linear connection we must first define the map $\sigma$
of (1.5) with ${\cal M} = \Omega^1(M_3^+)$.
Because of the identification (5.2) it can be considered as a
map from $M^+_3$ into itself and because of the bilinearity it is
necessarily of the form
$$
\sigma = \pmatrix{\mu &  0  & 0 \cr
                   0  & \mu & 0 \cr
                   0  &  0  & -1},
$$
where $\mu \in \bbbc$. The $-1$ in the lower right corner is imposed by
the condition (1.10).

It can be shown (Madore {\it et al.} 1995) that for each such $\sigma$
there is a unique linear connection given by the covariant derivative
(2.13). That is, necessarily $\tau \equiv 0$. Let $e$ be the unit in
$M_1$ considered as generator of $\Omega^2(M_3^+)$. The expression
$d\theta + \theta^2$ is given by
$$
d\theta + \theta^2 = e.
$$
Therefore we have
$$
{\rm Curv}_L (\xi) = e \otimes \xi.                               \eqno(5.3)
$$

To construct $J$ it is convenient to fix a vector-space basis for
$\Omega^1(M_3^+)$. We introduce the (unique) upper-triangular matrices
$\eta_1$ and $\eta_2$ such that $\theta = \eta_1 - \eta_1^*$ and such
that
$$
\eta_i \eta_j^* = 0, \qquad \eta_i^* \eta_j = \delta_{ij} e.
$$
We find (Madore {\it et al.} 1995) that
$$
\matrix{
\nabla^2 \eta_1 = 0,                                  \hfill
&\nabla^2 \eta_2 = 0,                                 \hfill\cr
\nabla^2 \eta^*_1 = -(\mu + 1) e \otimes \eta^*_1,    \hfill
&\nabla^2 \eta^*_2 = - e \otimes \eta^*_2.            \hfill
}                                                                 \eqno(5.4)
$$

Since there is an element $u$ of the algebra such that $\eta_2 = u \eta_1$
it is obvious that the map $\nabla^2$ is right-linear only in the
degenerate case $\mu = 0$. In this case $J = 0$ and
$$
{\rm Curv} = {\rm Curv}_L.                                        \eqno(5.5)
$$
Otherwise it is easy to see that
$$
J = \Omega^2(M_3^+) \otimes_{M^+_3} \Omega^1(M_3^+)
$$
and therefore that
$$
{\rm Curv} \equiv 0.                                              \eqno(5.6)
$$
It is difficult to appreciate the meaning of this result since there is
only one connection for each value of $\mu$. However (5.4) does not
appear to define a curvature which is any the less flat for generic
$\mu$ than for the special value $\mu = 0$. As we have defined it
${\rm Curv}$ does not perhaps contain enough information to characterize
a general noncommutative geometry.

{}From (5.4) one sees that for all values of $\mu$ there is a subalgebra
of $M_3^+$ with respect to which $\nabla^2$ is right-linear. It consists
of those elements which leave invariant the vector sub-spaces of
$\Omega_1(M_3^+)$ defined by $\eta_1$ and $\eta_2$. That is, it is the
algebra
$$
M_1 \times M_1 \times M_1 \subset M_3^+.
$$

As in Section~1 we define
$M_3^{+e} = M_3^+ \otimes_\bbbc M_3^{+ {\rm op}}$. A general element
$\xi$ of $\Omega^1(M_3^+)$ can be written in the form
$$
\xi = \pmatrix{   0    &      0     &  \xi_{13}  \cr
                  0    &      0     &  \xi_{23}  \cr
              \xi_{31} &  \xi_{32}  &      0       },             \eqno(5.7)
$$
where the $\xi_{ij}$ are arbitrary complex numbers. The map
$$
\xi \mapsto \pmatrix{0 & \xi_{13} & 0  \cr
                     0 & \xi_{23} & 0  \cr
                     0 &    0     & 0} \otimes
            \pmatrix{0 & 0 & 0 \cr
                     0 & 0 & 0 \cr
                     0 & 0 & 1} +
            \pmatrix{0 & 0 & 0 \cr
                     0 & 0 & 0 \cr
                     0 & 0 & 1} \otimes
            \pmatrix{  0     &    0     & 0 \cr
                    \xi_{31} & \xi_{32} & 0 \cr
                       0     &      0   & 0}                     \eqno(5.8)
$$
identifies $\Omega^1(M_3^+)$ as a sub-bimodule of the free
$M_3^{+e}$-module of rank 1 and the projector
$$
P = \pmatrix{0 & 0 & 0  \cr
             0 & 1 & 0  \cr
             0 & 0 & 0} \otimes
    \pmatrix{0 & 0 & 0  \cr
             0 & 0 & 0  \cr
             0 & 0 & 1} +
    \pmatrix{0 & 0 & 0  \cr
             0 & 0 & 0  \cr
             0 & 0 & 1} \otimes
    \pmatrix{0 & 0 & 0  \cr
             0 & 1 & 0  \cr
             0 & 0 & 0}                                            \eqno(5.9)
$$
projects $M_3^{+e}$ onto $\Omega^1(M_3^+)$. One immediately
sees that by multiplication of $P$ on the right and left by elements of
$M_3^+$ one obtains all elements of $\Omega^1(M_3^+)$. The construction
of Section~2 can be used to construct by projection a covariant
derivative (2.19). In the present case we find that
$$
P (\theta \otimes P) = 0
$$
as it must be since we have already noticed that in the present case
$\tau \equiv 0$.  The covariant derivative (2.19) is identical to
that given by (2.10).

\beginsection 6 Braided-commutative algebras

As an example of a braided-commutative differential calculus we consider
the quantum plane with its $SL_q(2)$-covariant differential calculus
$\Omega^*$. It has been found (Dubois-Violette {\it et al.} 1995) that
there is a unique 1-parameter family of linear connections given by
the covariant derivative
$$
D \xi = \mu^4 x^i (x\eta - q y \xi) \otimes (x\eta - q y \xi),     \eqno(6.1)
$$
with $\mu$ a complex number. The corresponding $\sigma$ is uniquely
defined in terms of the $R$-matrix.  There are other linear connections
if we extend the algebra to include additional elements $x^{-1}$ and
$y^{-1}$.  For example consider the construction of Section~2 based on
the formula (2.8).  For arbitrary complex number $c$ define for
$q \neq 1$
$$
\theta = {1 \over 1 - q^{-2}} (y^{-1} \eta + c x y^{-2} \xi).      \eqno(6.2)
$$
We have then
$$
\xi^i = dx^i = - [\theta, x^i].                                    \eqno(6.3)
$$
If $c = 0$ then the differential $d$ is given on the entire algebra of
forms as a graded commutator with $\theta$.  Using $\theta$ we can
define $D_L$ and $D_R$ by (2.10) and a covariant derivative by (2.13).
As in Section~4 the curvature of this covariant derivative
vanishes. In fact for arbitrary $c$ we have $d\theta + \theta^2 = 0$.
This construction can be used for any generalized permutation which
satisfies the condition (1.10). There are many such $\sigma$. For
example if $i$ is a bimodule injection of $\Omega^2$ into
$\Omega^1 \otimes_{\cal A} \Omega^1$ which satisfies the condition
$\pi \circ i = 1$ then a generalized permutation is given by the formula
$\sigma = 1 - 2 i \circ \pi$ (Mourad 1995).

A linear connection has also been constructed on $GL_q(n)$ (Georgelin
{\it et al.} 1995). The differential calculus is constructed using a
1-form $\theta$ and a linear connection is given by the Formula~(2.13).

The construction of a bilinear curvature based on the projection (3.2)
is not interesting in the general braided-commutative case.  In this
case the right-module structure of $\Omega^1({\cal A})$ is determined in
terms of its left-module structure even though the forms do not commute
with the algebra. The construction can be modified however using the
braiding. There is then a morphism $\rho$ of the algebra
such that the vector space

$$
J_\rho = \{ \sum_i \big(\nabla^2 (\xi_i f_i) -
                        \nabla^2 (\xi_i) \rho(f_i) \big)
         \; \vert \; \xi_i \in \Omega^1({\cal A}), f_i \in {\cal A} \}
$$
vanishes identically. The curvature ${\rm Curv}_\rho$ defined by the obvious
modification of (3.3) is therefore left linear and right $\rho$-linear:
$$
{\rm Curv}_\rho(f\xi) = f {\rm Curv}_\rho(\xi), \qquad
{\rm Curv}_\rho(\xi f) = {\rm Curv}_\rho(\xi) \rho(f).          \eqno(6.4)
$$
In general, for each automorphism $\rho$ of the algebra, a curvature
${\rm Curv}_\rho$ can be defined; it would be of interest however
only in the case when $J_\rho$ vanishes.

\beginsection 7 The problem of curvature invariants

In Section~4 we considered a geometry with a module of 1-forms which was
free of rank $n^2-1$ as a left (and right) module. We noticed also, in
(4.13) that it can be written as a direct summand in a free bimodule of
rank 1. In Section~5 we considered the projective bimodule structure of
$\Omega^1(M_3^+)$.  In this Section we shall examine the projective
structure of $\Omega^1(M_3^+)$ as a left (and right) module in order to
see to what extent it is possible to express the geometry of Section~5
in the language of Section~4. If it were possible to do this it would be
possible to define curvature invariants as in Section~4.

Consider the element (5.7) of $\Omega^1(M_3^+)$. The map
$$
\xi \mapsto \pmatrix{0 & 0 &    0      \cr
                     0 & 0 &    0      \cr
                     0 & 0 & \xi_{31}},\quad
            \pmatrix{0 & 0 &    0      \cr
                     0 & 0 &    0      \cr
                     0 & 0 & \xi_{32}},\quad
            \pmatrix{0 & \xi_{13} & 0  \cr
                     0 & \xi_{23} & 0  \cr
                     0 &    0     & 0}                            \eqno(7.5)
$$
identifies $\Omega^1(M_3^+)$ as a submodule of the free module
$$
{\cal M} \equiv (M_3^+)^3 = M_3^+ \oplus M_3^+ \oplus M_3^+
$$
of rank 3. This imbedding respects the left-module structure of $M_3^+$.
It respects also the right-module structure if we identify
$$
f \mapsto \rho(f) = \pmatrix{1 & 0 & 0  \cr
                             0 & 1 & 0  \cr
                             0 & 0 & 1} \otimes f.                 \eqno(7.6)
$$
That is, under the right action by $M_3^+$ the element $f$ it acts on
the row vector and not on the matrix entries.

Define the projectors
$$
P_1 = P_2 = \pmatrix{0 & 0 & 0  \cr
                     0 & 0 & 0  \cr
                     0 & 0 & 1},\qquad
P_3 = \pmatrix{0 & 0 & 0  \cr
               0 & 1 & 0  \cr
               0 & 0 & 0}                                          \eqno(7.7)
$$
in $M_3^+$ and the projector
$$
P = P_1 \otimes \pmatrix{1 & 0 & 0  \cr
                         0 & 0 & 0  \cr
                         0 & 0 & 0} +
    P_2 \otimes \pmatrix{0 & 0 & 0  \cr
                         0 & 1 & 0  \cr
                         0 & 0 & 0} +
    P_3 \otimes \pmatrix{0 & 0 & 0  \cr
                         0 & 0 & 0  \cr
                         0 & 0 & 1}                                \eqno(7.8)
$$
in $M_3(M_3^+)$. Then $\xi P = \xi$. Let $\alpha$ be a general element
of ${\cal M}$, that is, a triplet of elements of $M_3^+$ written as in
(7.5) as a row vector. Then $\alpha P \in \Omega^1(M_3^+)$ and all
elements $\xi \in \Omega^1(M_3^+)$ can be obtained in this way; the module
$\Omega^1(M_3^+)$ is a projective left $M_3^+$-module:
$$
\Omega^1(M_3^+) = {\cal M} P.                                     \eqno(7.9)
$$
This defines a projection
$$
{\cal M} \buildrel p \over \longrightarrow \Omega^1(M_3^+)        \eqno(7.10)
$$
which is a left inverse of the imbedding (7.5).

Let $\theta^r$ be the canonical basis of ${\cal M}$:
$$
\theta^1 = (1,0,0), \quad
\theta^2 = (0,1,0), \quad
\theta^3 = (0,0,1)
$$
where the unit is the unit in $M_3^+$. We use a notation here which
parallels that of Section~4 (with $n=2$). In general however, for
$f \in M_3^+$,
$$
f \theta^r \neq \theta^r \rho(f) \equiv \theta^s \big(\rho(f)\big)^r_s.
$$
This is an essential difference with the geometry of Section~4.

Define
$$
\theta^r_P = \theta^r P.
$$
By this we mean that $\theta^r_P$ is the image of the triplet $\theta^r$
under the projection (7.10) which we again identify as an element of
${\cal M}$ by (7.5). An extension $\tilde \sigma$ of $\sigma$ is a map
$$
\Omega^1(M_3^+) \otimes_{M_3^+} {\cal M}
\buildrel \tilde \sigma \over \longrightarrow
{\cal M} \otimes_{M_3^+} \Omega^1(M_3^+)
$$
given by the action $\tilde \sigma (\theta^r_P \otimes \theta^s)$. It is
clear that this $\tilde \sigma$ will not be a simple permutation as in
(4.1).  The covariant derivative on ${\cal M}$ will be defined by
$$
\tilde D \theta^r = - \omega^r{}_{st} \, \theta^s_P \otimes \theta^t
                                                                 \eqno(7.11)
$$
analogous to (4.2) but with here the $\omega^r{}_{st}$ arbitrary
elements of $M_3^+$.

Using the projection $p$ we can define in particular a covariant
derivative $\tilde D$ on ${\cal M}$ which coincides with the image of
$D$ on $\Omega^1(M_3^+)$ by the requirement that the diagram
$$

\matrix{
\Omega^1(M_3^+)& \buildrel D \over \longrightarrow
&\Omega^1(M_3^+) \otimes_{M_3^+} \Omega^1(M_3^+)      \cr
p \uparrow \phantom{p} && \downarrow                  \cr
{\cal M}& \buildrel \tilde D\over \longrightarrow
&\Omega^1(M_3^+) \otimes_{M_3^+} {\cal M}
}
                         \eqno(7.12)
$$
be commutative. The down arrow on the right is an injection defined by
(7.5).  The covariant derivative $\tilde D \theta^r$ is defined then by
$$
\tilde D \theta^r = D \theta^r_P.                                \eqno(7.13)
$$

\medskip

{\it Acknowledgment:}\ Part of this research was carried out while the
authors were visiting the Erwin Schr\"odinger Institut. They would
like to thank the acting director Peter Michor for his hospitality and
enlightening discussions.
\parskip 3pt plus 1pt
\parindent=0cm

\bigskip
\centerline{\bf References}
\medskip

Bresser K., M\"uller-Hoissen F., Dimakis A., Sitarz A. 1995,
{\it Noncommutative Geometry of Finite Groups}, Preprint GOET-TP 95/95.

Chamseddine A.H., Felder G., Fr\"ohlich J. 1993, {\it Gravity in
Non-Commutative Geometry}, Commun. Math. Phys. {\bf 155} 205.

Connes A. 1986, {\it Non-Commutative Differential Geometry}, Publications
of the Inst. des Hautes Etudes Scientifique. {\bf 62} 257.

--- 1994, {\it Noncommutative Geometry}, Academic Press.

Connes A., Lott J. 1990, {\it Particle Models and Noncommutative Geometry},
in `Recent Advances in Field Theory', Nucl. Phys. Proc. Suppl. {\bf B18} 29.

--- 1992, {\it The metric aspect of non-commutative geometry},
Proceedings of the 1991 Carg\`ese Summer School, Plenum Press.

Cuntz A., Quillen D. 1995, {\it Algebra extensions and nonsingularity},
J. Amer. Math. Soc. {\bf 8} 251.

Dabrowski L., Hajac P., Landi G., Siniscalco P. 1995, {\it Biconnections
on quantum spaces}, SISSA \& ICTP Preprint.

Dubois-Violette M. 1988, {\it D\'erivations et calcul diff\'erentiel
non-commutatif}, C. R. Acad. Sci. Paris {\bf 307} S\'erie I 403.

Dubois-Violette M., Kerner R., Madore J. 1989, {\it Gauge bosons in a
noncommutative geometry}, Phys. Lett. {\bf B217} 485; {\it Classical
bosons in a noncommutative geometry}, Class. Quant. Grav. {\bf 6} 1709.

--- 1990, {\it Noncommutative differential geometry of matrix algebras},
J. Math. Phys. {\bf 31} 316.

Dubois-Violette M., Madore J., Masson T., Mourad J. 1995, {\it Linear
Connections on the Quantum Plane}, Lett. Math. Phys. {\bf 35} 351.

Dubois-Violette M., Michor P. 1994, {\it D\'erivations et calcul
diff\'erentiel non-commutatif II}, C. R. Acad. Sci. Paris {\bf 319}
S\'erie I 927.

Dubois-Violette M., Michor P. 1995, {\it Connections on Central
Bimodules}, Preprint LPTHE Orsay 94/100.

Dubois-Violette M., Masson T. 1995, {\it On the First-Order Operators in
Bimodules}, Preprint LPTHE Orsay 95/56.

Georgelin Y., Madore J., Masson T., Mourad J. 1995,
{\it  On the non-commutative Riemannian geometry of $GL_{q}(n)$},
Preprint LPTHE Orsay 95/51.

Karoubi M. 1981, {\it Connections, courbures et classes
caract\'eristiques en $K$-th\'eorie alg\'ebrique.}, Current trends in
algebraic topology, Part I, London, Ont.

Klim\v c\'ik C., Pompo\v s A., Sou\v cek. V. 1994, {\it Grading of
Spinor Bundles and Gravitating Matter in Noncommutative Geometry},
Lett. Math. Phys.  {\bf 30}, 259.

Koszul J.L. 1960, {\it Lectures on Fibre Bundles and Differential Geometry},
Tata Institute of Fundamental Research, Bombay.

Landi G., Nguyen Ai Viet, Wali K.C. 1994, {\it Gravity and
electromagnetism in noncommutative geometry}, Phys. Lett. {\bf B326} 45.

Madore J., Masson T., Mourad J. 1995, {\it Linear Connections on
Matrix Geometries}, Class. Quant. Grav. {\bf 12} 1429.

Mourad. J. 1995, {\it Linear Connections in Non-Commutative Geometry},
Class. Quant. Grav. {\bf 12} 965.

Sitarz A. 1994, {\it Gravity from Noncommutative Geometry}, Class.
Quant. Grav. {\bf 11} 2127.

Woronowicz S. L. 1987, {\it Twisted $SU(2)$ Group. An example of a
Non-Commutative Differential Calculus}, Publ. RIMS, Kyoto Univ. {\bf 23} 117.

\bye